\documentclass[11pt]{amsart}

\usepackage{latexsym,amssymb,amsmath,graphics}

\newcounter{rev}
\usepackage[dvips]{graphicx}

\def\cc{{\mathcal C}}

\def\ff{{\mathcal F}}

\def\ss{{\mathcal S}}

\def\ffi{\varphi}
\def\eps{\varepsilon}
\def\dst{\displaystyle}

\def\C{{\mathbb{C}}}

\def\N{{\mathbb{N}}}

\def\R{{\mathbb{R}}}

\newcommand{\scal}[1]{{\left\langle{#1}\right\rangle}}

\newenvironment{remark}[1][]{\vskip1pt\noindent\rm\textit{Remark.}\ }{\rm\vskip1pt}

\newenvironment{example}[1][]{\vskip1pt\noindent\rm\textit{Example.}\ }{\rm\vskip1pt}

\newtheorem{fact}{Fact}

\newtheorem{lemma}{Lemma}[section]
\newtheorem{proposition}[lemma]{Proposition}
\newtheorem{theorem}[lemma]{Theorem}
\newtheorem{corollary}[lemma]{Corollary}

\begin{document}

\title[Positive positive definite functions]{On the extremal rays of the cone of positive, positive definite functions}
\author{Philippe Jaming, Mat\'e Matolcsi \& Szil\'ard Gy. R\'ev\'esz}

\address{P.J.: Universit\'e d'Orl\'eans\\
Laboratoire MAPMO, CNRS, UMR 6628\\
Fédération Denis Poisson, FR 2964\\
B\^atiment de Mathématiques\\
BP 6759\\
45067 Orleans cedex 2\\
France}
\email{Philippe.Jaming@univ-orleans.fr}

\address{M.M.: R\'enyi Institute of Mathematics, 1053, Realtanoda u. 13-15, Budapest, Hungary}
\email{matomate@renyi.hu}
\address{Sz.R.: R\'enyi Institute of Mathematics, 1053, Realtanoda u. 13-15, Budapest, Hungary}
\email{revesz@renyi.hu}

\begin{abstract}
The aim of this paper is to investigate the cone of non-negative,
radial, positive-definite functions in the set of continuous
functions on $\R^d$.
Elements of this cone admit a Choquet integral representation in
terms of the extremals. The main feature of this article is to characterize some large classes
of such extremals. In particular, we show that there many other extremals than the gaussians,
thus disproving a conjecture of G. Choquet and that no reasonable conjecture can be made
on the full set of extremals.

The last feature of this article is to show that many characterizations of positive definite
functions available in the literature are actually particular cases of the Choquet integral representations
we obtain.
\end{abstract}


\thanks{The authors wish to thank E.\,G.\,F. Thomas for valuable conversations and 
G. Godefroy for mentioning the problem to them.}
\maketitle

\tableofcontents

\section{Introduction}

Positive definite functions appear in many areas of mathematics, ranging from number theory
to statistical applications. Since Bochner's work,
these functions are known to be characterized as having a non-negative Fourier transform.

Before going on, let us first fix some notations. We will define the Fourier
transform of a function $f\in L^1(\R^d)$ by
$$
\ff_d f(\xi)=\widehat{f}(\xi)=\int_{\R^d}f(x)e^{2i\pi\scal{x,\xi}}\,\mbox{d}x
$$
and extend this definition both to bounded measures on $\R^d$ and
to $L^2(\R^d)$ in the usual way. Here $\scal{\cdot,\cdot}$ is the
scalar product on $\R^d$ and $|\cdot|$ the Euclidean norm.

A continuous function $f$ is said to be \emph{positive definite} if,
for every integer $n$, for all $x_1,\ldots,x_n\in\R^d$,
the $n\times n$ matrix $[f(x_j-x_k)]_{1\leq j,k\leq n}$ is positive definite,
that is, if
\begin{equation}
\label{def.posdef}
\sum_{j,k=1}^nc_j\overline{c_k}f(x_j-x_k)\geq 0\qquad\mbox{for all }c_1,\ldots,c_n\in\C.
\end{equation}
Then Bochner's Theorem \cite{Boc} shows that $f$ is positive
definite if and only if $f=\widehat{\mu}$ for some positive
bounded Radon measure on $\R^d$ (a probability measure if we
further impose $f(0)=1$). There are many proofs of Bochner's
theorem, the nearest to the subject of this paper being based on
the Choquet Representation Theorem, due to Bucy and Maltese
\cite{BM}, {\it see also} \cite{Be,Ch,Ph}. Let us recall the main
features of these proofs, and thereby also some definitions.
An element $f$ of a cone $\Omega\ni f$ is an extremal ray generator of $\Omega$
(or simply an extremal) if $f=f_1+f_2$ with $f_1,f_2\in\Omega$
implies $f_1=\lambda f$, $f_2=\mu f$, $\lambda,\mu\geq 0$. The
first step in the proof of Bochner's theorem is then to show that
the characters $e^{2i\pi\scal{x,\xi}}$ are the only extremal rays
of the cone of positive definite functions. The second step is to
show that the cone of positive definite functions on $\R^d$
is well capped, {\it i.e.} is the union of caps
(compact, convex subsets $C$ of $\Omega$ such that
$\Omega\setminus C$ is still convex). It then follows from the
work of Choquet that every element of such a cone is an integral
over extreme points with respect to a conical measure. For more
details, we refer to the references given previously.

\medskip

Let us now note that Bochner's Theorem, though being powerful when
one wants to construct positive definite functions, may be
difficult to use in practice. This mainly comes from the fact that
explicit computations of Fourier transforms are generally
impossible. For instance, it is not known precisely for which
values of $\lambda$ and $\kappa$ the function
$(1-|x|^\lambda)_+^\kappa$ is positive definite on $\R^d$. This
problem is known as the Kuttner-Golubov problem and we refer to
\cite{Gnproc} for more details and the best known results to date.
To overcome this difficulty, one seeks concrete and easily
checkable criteria that guarantee that a function is positive
definite. The most famous such criterion is due to P\'olya which
shows that a bounded continuous function on $\R$ which is convex
on $[0,+ \infty)$, is positive definite. More evolved criteria may
be found in the literature ({\it see} Section \ref{sec:criteria}
for more details).

As it turns out, the functions so characterized are not only
positive definite but also non-negative. We will call such
functions \emph{positive positive definite}. Such functions
appear in many contexts. To give a few examples where the reader
may find further references, let us mention various fields such as
approximation theory \cite{Bu}, spatial statistics \cite{Bon},
geometry of Banach spaces \cite{Ko}, and physics \cite{GS}.
Despite a call to study such functions by P. L\'evy in \cite{Le}
they seem not to have attracted much attention so far. To our knowldge,
there are only two papers specifically devoted to the subject in the litterature: \cite{GiPe}
which is motivated by applications in physics and the
(unpublished) paper \cite{Bor} which is motivated by problems in number
theory.

\medskip

Before going on with the  description of the main features of this
paper, we will need some further notations. For $d\geq 1$ we
define $\cc_r(\R^d)$ to be the space of radial continuous
functions, and we stress the fact that in the sequel we will only
consider \emph{radial} functions in dimensions higher than 1. Now,
let
$$
\Omega_{d}=\{f\in \cc_r(\R^d)\,: f\mbox{ is }\mbox{positive
definite}\}
$$
and $\Omega_{d}^+=\{f\in\Omega_{d}\,: f\geq0\}$.
Note that in dimension $1$, a positive positive definite function is even, so there
is no restriction when considering radial functions in this case.

Next, $\Omega_{d}^+$ is a closed convex sub-cone of $\Omega_d$ in $\cc_r(\R^d)$
(endowed with the weak $*$ $\sigma(L^1,L^\infty)$ topology). For
$f\in\Omega_{d}^+$, we denote by
$I(f)=\Omega_d^+\cap(f-\Omega_d^+)$ the \emph{interval} generated
by $f$. Then $f$ is an \emph{extremal ray
generator} if $I(f)=\{\lambda f\,:\
0\leq\lambda\leq 1\}$. As the cone $\Omega_d$ is well-capped, so
is $\Omega_d^+$. Therefore,
Choquet Theory applies and every positive positive definite
function admits an integral representation over extremals. It is
therefore a natural task to determine the extremals of
$\Omega_d^+$.

We are unfortunately unable to fulfill this task
completely. One difficulty is that among the extremals of the cone of positive
definite functions, only the trivial character $1$ is still in the cone of
positive positive definite functions and is of course an extremal of it.
We are nevertheless  able to describe large classes of
extremals. One such class is included in the compactly supported
extremals of $\Omega_{d}^+$. In this case, we show that if an
element of $\Omega_{d}^+$ with compact support has a Fourier
transform whose holomorphic extension to $\C$ has only real
zeroes, then this element is an extremal ray generator. One would
be tempted to conjecture that this describes all compactly
supported extremals, but we show that this is not the case.
Nevertheless, this theorem allows to show that many examples of
functions used in practice are extremals, as for instance the
functions $(1-|x|^2)_+^\alpha*(1-|x|^2)_+^\alpha$ for suitable
$\alpha$'s.

The next class of functions we investigate is that of Hermite
functions, that is, functions of the form $P(x)e^{-\lambda x^2}$,
$P$ a polynomial. This is a natural class to investigate since the
elements of the intervals they generate consist of functions of
the same form. This will be shown as a simple consequence of
Hardy's Uncertainty Principle. Further, a conjecture that oral
tradition \cite{Godefroy} attributes to G. Choquet (although P. L\'evy may be
another reasonable source of the conjecture) states that the only
extremals are the gaussians. This conjecture is false, as our results on
compactly supported extremals show. We will construct more counter-examples
by describing precisely the positive positive definite
functions of the form $P(x)e^{-\pi x^2}$ where $P$ is a polynomial
of degree $4$ and showing that this class contains extremal ray
generators. This further allows us to construct extremals of the
form $P(x)e^{-\pi x^2}$ with $P$ polynomials of arbitrary high
degree.

\medskip

Finally, we show that most (sufficient) characterizations of positive definite functions
actually characterize positive positive definite functions and are actually particular cases of
Choquet representations. More precisely, it is easy to see that if $\ffi$ is an extremal ray
generator in $\Omega_{d}^+$ then so is $\ffi_t(x)=\ffi(tx)$. It follows that
\begin{equation}
\label{eq:mixing}
\dst\int_0^{+\infty}\ffi_t\,\mbox{d}\mu(t)
\end{equation}
is a positive positive definite function (for suitable $\mu$) that
is obtained by a Choquet representation with a measure supported
on the family of extremals $\{\ffi_t\}$. For instance:

--- as a particular case of our theorem concerning compactly supported functions, we obtain that the
function $\ffi(x)=(1-|x|)_+$ is extremal. P\'olya's criterium may
be seen as a characterization of those functions which may be
written in the form \eqref{eq:mixing}.

--- a criterium for deciding which functions may be written in the form \eqref{eq:mixing}
with $\ffi(x)=(1-|x|^2)_+*(1-|x|^2)_+$ has been obtained by Gneiting.

\medskip

The article is organized as follows. In the next section, we
gather preliminaries on positive definite functions. We then turn
to the case of compactly supported functions in Section
\ref{sec:compsupp}, followed by the case of Hermite functions. We
then devote Section \ref{sec:criteria} to link our results with
various criteria available in the literature.  We conclude the
paper with some open questions.

\section{Preliminaries}

\subsection{Basic facts}


\begin{fact}[Invariance by scaling]\ \\
Let $f$ be a continuous function and $\lambda>0$. We write
$f_\lambda(x)=f(\lambda x)$. Then $f\in \Omega_d^+$ if and only if
$f_\lambda\in\Omega_d^+$. Moreover $f$ is an extremal ray
generator if and only if $f_\lambda$ is extremal.
\end{fact}

Another consequence is the following lemma:

\begin{lemma}
\label{lem:extcrit} Let $\omega\in\Omega_d^+$ and $\nu$ be a
positive bounded measure on $(0,+\infty)$. Define
$$
F(x)=\int_0^{+\infty}\omega(x/t)\,\mbox{d}\nu(t).
$$
Then $F\in\Omega_d^+$. Moreover, write $\omega(x)=\omega_0(|x|)$
and assume that $\omega_0$ is non-increasing on $(0,+\infty)$ and
(strictly) decreasing in a neighborhood of $0$. Then $F$ is an
extremal ray generator if and only if $\omega$ is an extremal ray
generator and $\nu=\delta_a$ is a Dirac mass for some
$a\in(0,+\infty)$.
\end{lemma}

\begin{proof}
Note that the fact that $F$ is continuous follows from Lebesgue's Theorem since $\omega$ is bounded
and continuous. It is then also obvious that $F$ is positive and positive definite.

If $\omega$ is an extremal ray generator and $\nu=\delta_a$ then $F(x)=\omega(x/a)$ is clearly an extremal
ray generator. It is also immediate that if $\omega$ is not extremal then $F$ is not extremal.

Finally, assume that the support of $\nu$ contains at least 2
points $a<b=a+3\eta$ and let $\psi$ be a continuous non-increasing function such that
$0\leq\psi\leq 1$ and $\psi(x)=1$ on $(0,a+\eta)$ while $\psi=0$ on $(a+2\eta,+\infty)$. Define
$$
F_1(x)=\int_0^{+\infty}\omega(x/t)\psi(t)\,\mbox{d}\nu(t)
\quad\mbox{ and }\quad
F_2(x)=\int_0^{+\infty}\omega(x/t)\bigl(1-\psi(t)\bigr)\,\mbox{d}\nu(t).
$$
The measures $\psi(t)\,\mbox{d}\nu(t)$ and $\bigl(1-\psi(t)\bigr)\,\mbox{d}\nu(t)$
are bounded, so that $F_1,F_2$ are continuous, $F_1$ and
$F_2$ are positive positive definite, $F=F_1+F_2$.

Assume now towards a contradiction that $F_1=\lambda F$ with
$0\leq\lambda\leq 1$. From the assumption on $\psi$, we easily
deduce that $0<\lambda<1$. Let $t_0$ be a solution of
$\psi(t)=\lambda$. From $\dst\int_{0}^{+\infty}
\omega(x/t)\bigl(\psi(t)-\lambda\bigr)\,\mbox{d}\nu(t)=0$ we get
that
\begin{equation}
\label{eq:psi}
\int_0^{t_0}\omega(x/t)\bigl(\psi(t)-\lambda\bigr)\,\mbox{d}\nu(t)=
\int_{t_0}^{+\infty}\omega(x/t)\bigl(\lambda-\psi(t)\bigr)\,\mbox{d}\nu(t).
\end{equation}
In particular, for $x=0$, we obtain
$$
\int_0^{t_0}\bigl(\psi(t)-\lambda\bigr)\,\mbox{d}\nu(t)
=\int_{t_0}^{+\infty}\bigl(\lambda-\psi(t)\bigr)\,\mbox{d}\nu(t).
$$
Now, $\chi_{(0,t_0)}\bigl(\psi(t)-\lambda\bigr)\,\mbox{d}\nu(t)$
and
$\chi_{(t_0,+\infty)}\bigl(\lambda-\psi(t)\bigr)\,\mbox{d}\nu(t)$
are positive non-zero measures. Then, for $|x_0|>0$  small enough to
have that $\omega$ is strictly decreasing on $(0,|x_0|/t_0)$, we get
\begin{eqnarray*}
\int_0^{t_0}\omega(x_0/t)\bigl(\psi(t)-\lambda\bigr)\,\mbox{d}\nu(t)
&=&\int_0^{t_0}\omega_0(|x_0|/t)\bigl(\psi(t)-\lambda\bigr)\,\mbox{d}\nu(t)\\
&\leq&\omega_0(|x_0|/t_0)\int_0^{t_0}\bigl(\psi(t)-\lambda\bigr)\,\mbox{d}\nu(t)\\
&=&\omega_0(|x_0|/t_0)\int_{t_0}^{+\infty}\bigl(\lambda-\psi(t)\bigr)\,\mbox{d}\nu(t)\\
&<&\int_{t_0}^{+\infty}\omega_0(|x_0|/t)\bigl(\lambda-\psi(t)\bigr)\,\mbox{d}\nu(t)\\
&=&\int_{t_0}^{+\infty}\omega(x_0/t)\bigl(\lambda-\psi(t)\bigr)\,\mbox{d}\nu(t)
\end{eqnarray*}
a contradiction, so that $F_1\ne \lambda F$ and $F$ is not an extremal ray generator.
\end{proof}

\begin{fact}[Invariance under products and convolution]\ \\
If $f,g\in\Omega_d^+$ then $fg\in\Omega_d^+$.
Further, if $f,g$ are also in $L^2$ (say) then $f*g\in\Omega_d^+$.

Further if either $f$ or $g$ (resp. $\widehat{f}$ or $\widehat{g}$) does not vanish, then
, for $fg$ (resp. $f*g$) to be extremal, it is necessary that
both $f$ and $g$ are extremal.
\end{fact}

Indeed, assume that $g$ is not extremal and write
$g=g_1+g_2$ with $g_1/g,g_2/g$ not constant, then $fg=fg_1+fg_2$. Now, if $fg_1=\lambda fg$
then $g_1=\lambda g$ on the support of $f$.

The converse is unclear and probably false. A possible counter-example may be constructed
as follows. Assume there is a compactly supported extremal $f$ such that $f^2$ is also
extremal. Without loss of generality, we may assume that
$f$ is supported in $[-1,1]$. Let $g=(4\delta_0+\delta_{-4\pi}+\delta_{-4\pi})*f$
then $fg=4f^2$ would be extremal but $g$ is not.

\medskip

\subsection{Bessel functions and Fourier transforms}\ \\
Results in this section can be found in most books on Fourier
analysis, for instance \cite[Appendix B]{Gr}.

Let $\lambda$ be a real number with $\lambda>-1/2$. We define the
Bessel function $J_\lambda$ of order $\lambda$ on $(0,+\infty)$ by
its \emph{Poisson representation formula}
$$
J_\lambda(x)=\frac{x^\lambda}{2^\lambda
\Gamma\left(\lambda+\frac{1}{2}\right)\Gamma\left(\frac{1}{2}\right)}
\int_{-1}^1 e^{isx}(1-s^2)^\lambda\frac{\mbox{d}s}{\sqrt{1-s^2}}.
$$
Let us define $\mathcal{J}_{-1/2}(x)=\cos x$ and for
$\lambda>-1/2$,
$\mathcal{J}_\lambda(x):=\frac{J_\lambda(x)}{x^\lambda}$. Then
$\mathcal{J}_\lambda$ extends to an even entire function of order
$1$ and satisfies $\mathcal{J}_\lambda(x)$ real and $\mathcal{J}_\lambda(ix)>0$
for all $x\in\R$. It is also known that $\mathcal{J}_\lambda$ has
only real simple zeroes.

As is well-known, if $f$ is a radial function given by
$f(x)=f_0(|x|)$, then its Fourier transform is given by
$$
\widehat{f}(\xi)=\mathcal{J}_{\frac{d}{2}-1}(f_0)(|\xi|)
$$
where
$$
\mathcal{J}_\lambda(f_0)(t)=(2\pi)^{\lambda+1}\int_0^{+\infty}f_0(r)
\mathcal{J}_\lambda(2\pi r t)r^{2\lambda+1}\,\mbox{d}r.
$$
Bochner's Theorem has been extended to radial continuous positive
definite functions by Schoenberg \cite[page 816]{Sc} ({\it see
also} \cite{SvP})\,: a function $\ffi$ is radial positive definite
and radial on $\R^d$, $d\geq 2$ if and only if there exists a
positive bounded measure $\mu$ on $(0,+\infty)$ such that
$\ffi(x)=\ffi_0(|x|)$ with
$$
\ffi_0(r)=r^{-\frac{d}{2}+1}\int_0^{+\infty}J_{\frac{d}{2}-1}(2\pi
rs)\,\mbox{d}\mu(s)
=(2\pi)^{\frac{d}{2}-1}\int_0^{+\infty}\mathcal{J}_{\frac{d}{2}-1}(2\pi
rs)\,s^{\frac{d}{2}-1}\mbox{d}\mu(s).
$$
For $d=1$ this coincides with Bochner's theorem.

We will also use the following well-known result: denote $|t|_+=t$
or $0$ according to $t>0$ or not. Let
$m_\alpha(x)=(1-|x|^2)^\alpha_+$. Then
$$
\widehat{m_\alpha}(\xi)=\frac{\Gamma(\alpha+1)}{\pi^\alpha}
\frac{J_{\frac{d}{2}+\alpha}(2\pi|\xi|)}{|\xi|^{\frac{d}{2}+\alpha}}
=2^{\frac{d}{2}+\alpha}\pi^{\frac{d}{2}}\Gamma(\alpha+1)\mathcal{J}_{\frac{d}{2}+\alpha}(2\pi|\xi|).
$$

\section{Compactly supported positive positive definite functions}
\label{sec:compsupp}

In this section, we consider compactly supported positive positive
definite functions. It is natural to look for extremals inside
this class of functions because of the following (trivial)
lemma\,:

\begin{lemma}\label{lem:intcomp}
Let $f$ be a continuous radial positive positive definite function
with compact support. Then the interval $I(f)$ contains only
positive positive definite functions with support included in
$\mathrm{supp}\,f$.

Moreover, if we write $B(0,a)$ for the smallest ball containing
$\mathrm{supp}\,f$, {\it i.e.}
$\mathrm{conv}\,\mathrm{supp}\,f=B(0,a)$, and if $f=g+h$ with
$g,h\subset I(f)$, then at least one of
$\mathrm{conv}\,\mathrm{supp}\,g$ and
$\mathrm{conv}\,\mathrm{supp}\,h$ is $B(0,a)$.
\end{lemma}

\begin{proof}
First if $g\in I(f)$ then $0\leq g\leq f$ so that
$\mbox{supp}\,g\subset\mbox{supp}\,f$.

Next, as $f$ is radial, there exists $a$ such that
$\mbox{conv}\,\mbox{supp}\,f=B(0,a)$. Assume now that $f=g+h$ with
$g,h\in I(f)$. As $g$ (resp. $h$) is radial, the convex hull of its support is a
ball and we denote it by $B(0,b)$ (resp. $B(0,c)$). As $g$ and $h$
are both non-negative, the convex hull of the support of $g+h=f$ is then 
$B\bigl(0,\max(b,c)\bigr)=B(0,a)$ thus the claim.
\end{proof}

We may now prove the following:

\begin{theorem}\label{th.compextray}
Let $f\in\cc(\R^d)$ be a compactly
supported positive positive definite radial function. Write
$f(x)=f_0(|x|)$ where $f_0$ is a compactly supported function on
$\R_+$ so that $\mathcal{J}_{\frac{d}{2}-1}(f_0)$ extends
analytically to $\C$.

Assume that $\mathcal{J}_{\frac{d}{2}-1}(f_0)$ has only real
zeroes, then $f$ is an extreme ray generator in the cone of continuous
positive positive definite radial functions.

Moreover, assume that $\mathcal{J}_{\frac{d}{2}-1}(f_0)$ has only
a finite number $N$ of non-real zeroes, and let $g\in I(f)$. If we
write $g=g_0(|x|)$, then $\mathcal{J}_{\frac{d}{2}-1}(g_0)$ has at
most $N$ non-real zeroes.
\end{theorem}

\begin{proof}
Let us write $f=g+h$ with $g,h$ radial positive positive definite functions.
Then $g$ and $h$ are also compactly supported.

Write $g(x)=g_0(|x|)$ and $h(x)=h_0(|x|)$. It follows that
$\mathcal{J}_{\frac{d}{2}-1}(f_0)$,
$\mathcal{J}_{\frac{d}{2}-1}(g_0)$ and
$\mathcal{J}_{\frac{d}{2}-1}(h_0)$ all extend to entire functions
of order $1$. From Hadamard's factorization theorem, we may write
$\mathcal{J}_{\frac{d}{2}-1}(f_0)$ as
\begin{equation}
\label{eq:repf}
\mathcal{J}_{\frac{d}{2}-1}(f_0)(z)=z^ke^{az+b}
\prod_{\zeta\in \mathcal{Z}(f)}\left(1-\frac{z}{\zeta}\right)\exp\frac{z}{\zeta}
\end{equation}
where $\mathcal{Z}(f)$ is the set of non-zero zeroes of $\mathcal{J}_{\frac{d}{2}-1}(f_0)$.
Let us further note that
\begin{enumerate}
\renewcommand{\theenumi}{\roman{enumi}}
\item $\mathcal{J}_{\frac{d}{2}-1}(f_0)$ is real if $z\in \R$, thus $a,b\in\R$ and if $\zeta\in\mathcal{Z}(f)$ then $\overline{\zeta}
\in\mathcal{Z}(f)$;

\item $\mathcal{J}_{\frac{d}{2}-1}(f_0)$ is non-negative if $z\in \R$, thus real zeroes
are of even order;

\item $\mathcal{J}_{\frac{d}{2}-1}(f_0)(0)=\widehat{f}(0)=\dst\int f\not=0$ since $f\geq 0$ thus $k=0$;

\item $\mathcal{J}_{\frac{d}{2}-1}(f_0)$ is even thus if $\zeta\in \mathcal{Z}(f)$ then
$-\zeta\in \mathcal{Z}(f)$ and $a=0$.
\end{enumerate}
We may thus simplify \eqref{eq:repf} to
$\mathcal{J}_{\frac{d}{2}-1}(f_0)(z)= \widehat{f}(0)E_f(z)P_f(z)$
with
$$
E_f(z)=\prod_{\zeta\in
\mathcal{Z}_+(f)}\left(1-\frac{z^2}{\zeta^2}\right)^2
$$
where $\mathcal{Z}_+(f)=\mathcal{Z}(f)\cap(0,+\infty)$
and
$$
P_f(z)=\prod_{\zeta\in
\mathcal{Z}_Q(f)}\left(1-\frac{z^2}{\zeta^2}\right)\left(1-\frac{z^2}{\overline{\zeta}^2}\right)
$$
where $\mathcal{Z}_Q(f)=\{\zeta\in\mathcal{Z}(f)\,:\
\mbox{Re}\,\zeta>0\ \&\ \mbox{Im}\,\zeta>0\}$. Similar expressions
hold for $\mathcal{J}_{\frac{d}{2}-1}(g_0)$ and
$\mathcal{J}_{\frac{d}{2}-1}(h_0)$. It should also be noticed that
both $E_f$ and $P_f$ are non-negative on the real and the
imaginary axes.

Now, let us assume that $\mathcal{J}_{\frac{d}{2}-1}(f_0)$ has only finitely many non-real zeroes,
so that $P_f$ is a polynomial. In particular, there exists an integer $N$ and
a constant $C$ such that $|P_f(z)|\leq C(1+|z|)^N$.

As
$\mathcal{J}_{\frac{d}{2}-1}(f_0)=\mathcal{J}_{\frac{d}{2}-1}(g_0)+\mathcal{J}_{\frac{d}{2}-1}(h_0)$
with $\mathcal{J}_{\frac{d}{2}-1}(h_0)\geq 0$, we get $0\leq
\mathcal{J}_{\frac{d}{2}-1}(g_0)(z)\leq\mathcal{J}_{\frac{d}{2}-1}(f_0)(z)$
for $z$ \emph{real}. It follows that
$\mathcal{Z}_+(f)\subset\mathcal{Z}_+(g)$, with multiplicity.
Thus, we may partition the multiset
$\mathcal{Z}_+(g)=\mathcal{Z}_+(f)\cup \mathcal{Z}'(g)$.

From Hadamard's factorization, it follows that
$$
\mathcal{J}_{\frac{d}{2}-1}(g_0)(z)={\widehat{g}(0)}E_f(z)
\prod_{\zeta\in
\mathcal{Z}'(g)}\left(1-\frac{z^2}{\zeta^2}\right)\cdot P_g(z).
$$
So, we have written
$\mathcal{J}_{\frac{d}{2}-1}(g_0)(z)=G(z)E_f(z)$ where $G$ is an
entire function of order at most $1$. Note that $0\leq
\mathcal{J}_{\frac{d}{2}-1}(g_0)(x)\leq\mathcal{J}_{\frac{d}{2}-1}(f_0)(x)$,
so that $0\leq |G(x)|\leq C(1+|x|)^N$ for $x$ real.

Further, from the positivity of $f$ and $g$,
\begin{equation}\label{imag}
0< \mathcal{J}_{\frac{d}{2}-1}(g_0)(it)=\int_{B(0,a)}g(x)e^{2\pi
tx}\,\mbox{d}x\leq \int_{B(0,a)} f(x)e^{2\pi
tx}\,\mbox{d}x=\mathcal{J}_{\frac{d}{2}-1}(f_0)(it).
\end{equation}
 It follows that $|G(z)|$ is also bounded by
$C(1+|z|)^N$ on the imaginary axis. By Phragm\'en-Lindel\"of's
Principle $G$ is bounded by $2^NC(1+|z|)^N$ over each of the four
quadrants $Q_{\eps_1,\eps_2}=\{\eps_1\mbox{Re}\,z\geq0\ \&\
\eps_2\mbox{Im}\,z\geq0\}$, $\eps_1=\pm1,\eps_2=\pm1$. From
Liouville's Theorem, we thus get that $G$ is a polynomial of
degree at most $N$.

In particular, if $N=0$, then $G$ is a constant and
$\mathcal{J}_{\frac{d}{2}-1}(g_0)=\lambda \mathcal{J}_{\frac{d}{2}-1}(f_0)$
thus $g=\lambda f$. It follows that $h=(1-\lambda)f$, and $f$ is an extremal ray generator.
\end{proof}

\begin{remark}\ \\
Note that in the course of the proof we have shown that if the
Fourier transform $\widehat{f}$ of a compactly supported positive
positive definite function $f$ has a non real zero $\zeta$, then
$\zeta$ is not purely imaginary (see \eqref{imag}) and $-\zeta$,
$\overline{\zeta}$, $-\overline{\zeta}$ are also zeroes of
$\widehat{f}$.
\end{remark}

It should be noted that checking whether a particular function is an extremal positive positive definite
function may be difficult in practice. Neverteless, we will now give a few examples.

\begin{example}\ \\
Let us consider the characteristic (indicator) function
$\chi_{[-1,1]}$ of the interval $[-1,1]$ in $\R$.
Then $p(x):=\chi_{[-1,1]}*\chi_{[-1,1]}$ is an extremal ray generator
since its Fourier transform is $\widehat{p}(\xi)=\left(\frac{\sin 2\pi\xi}{\pi\xi}\right)^2$.

The celebrated positive definiteness criteria of P\'olya characterizes those functions
that may be written in the form $\int p(x/r)\,\mbox{d}\mu(r)$ with $\mu$ a positive measure
({\it see} Section \ref{sec:polya}).

Further examples are obtained by convolving $p$'s:
$$
\chi_{[-r_1/2,r_1/2]}*\chi_{[-r_1/2,r_1/2]}*\cdots*\chi_{[-r_{k}/2,r_{k}/2]}*\chi_{[-r_{k}/2,r_{k}/2]}.
$$
\end{example}

Another class of examples is given by the following:

\begin{corollary}
For $\alpha>-1/2$ define the function $m_\alpha$ on $\R^d$ by
$$
m_\alpha(x)=(1-|x|^2)_+^\alpha:=\begin{cases}(1-|x|^2)^\alpha&\mbox{if }|x|<1\\
0&\mbox{otherwise}\end{cases}.
$$
Then, $m_\alpha*m_\alpha$ is a continuous positive positive definite
function that is an extremal ray generator.
\end{corollary}

This covers the previous example and extends it to higher dimension since $m_0$ is the characteristic (indicator)
function of the unit ball of $\R^d$.

\medskip

\begin{remark}
If we were considering positive definite tempered distibutions, then this result stays true
for $\alpha>-1$.
\end{remark}

\begin{proof}
First of all, note that if $\alpha\geq 0$, $m_\alpha\in L^\infty$
and for $-1<\alpha<0$, $m_\alpha\in L^p$ for all $\alpha>-1/p$. It
follows that $m_\alpha*m_\alpha$ is well defined, has compact
support and is continuous for all $\alpha>-1/2$.

Further, as is well-known ({\it see e.g.} \cite[Appendix B5]{Gr}):
$$
\widehat{m_\alpha}(\xi)=
\frac{\Gamma(\alpha+1)}{\pi^\alpha}\frac{J_{d/2+\alpha}(2\pi|\xi|)}{|\xi|^{d/2+\alpha}}.
$$
It follows that the Fourier transform of $m_\alpha*m_\alpha$ is given by
$$
\widehat{m_\alpha*m_\alpha}(\xi)=\frac{\Gamma(\alpha+1)^2}{\pi^{2\alpha}}\frac{J_{d/2+\alpha}(2\pi|\xi|)^2}{|\xi|^{d+2\alpha}}.
$$
But, from a Theorem of Hurwitz \cite{Hu} ({\it see} \cite[{\bf
15.27} page 483]{Wa}), $J_{d/2+\alpha}$ has only positive zeroes
since $d/2+\alpha>-1$. The corollary thus follows from the
previous theorem.
\end{proof}

\begin{example}\ \\
A particular case of the previous result is that the function
$w(x):=m_1(x)*m_1(x)=(1-x^2)_+*(1-x^2)_+$ on $\R$ is an extremal
ray generator. The function $w$ has been introduced in the study
of positive definite functions by Wu \cite{Wu} in the context of
radial basis function interpolation. A simple criteria for writing
a function in the form $\int w(x/r)\,\mbox{d}\mu(r)$ with $\mu$ a
positive bounded measure has been recently obtained by Gneiting
\cite{Gn1} ({\it see} Section \ref{sec:gneiting}).

Further examples are then obtained by taking scales $w_r$ of $w$ and convolving several $w_r$'s together:
$$
(1-x^2/r_1^2)_+*(1-x^2/r_1^2)_+*\cdots*(1-x^2/r_k^2)_+*(1-x^2/r_k^2)_+
$$
\end{example}

\begin{remark}\\
When $\lambda\to \infty$ we have
$m_\lambda(\sqrt{\pi}x/\sqrt{\lambda})\to e^{-\pi x^2}$ pointwise,
in $L^2$, and uniformly on compact sets.
\end{remark}

It would be tempting to conjecture that every compactly supported
extreme ray generator is covered by Theorem \ref{th.compextray},
i.e. it is of the form that its Fourier transform has only real
zeroes. Nevertheless, this is not the case:

\begin{proposition}
There exists a continuous compactly supported extreme ray
generator $f$ of $\Omega_1^+$ such that $\widehat{f}$ has non-real
zeroes.
\end{proposition}

\begin{proof}
Consider $\ffi(x)=(1-x^2)_+^2*(1-x^2)_+^2$. A cumbersome
computation (by computer) shows that
$$
\ffi(x)=\begin{cases}\frac{1}{630}(2-x)^5(x^4+10|x|^3+36x^2+40|x|+16)&\mbox{if }|x|<2\\
0&\mbox{otherwise}\end{cases}.
$$
Let $r>0$ and $0<\theta<\pi/2$ and $\zeta=\frac{r}{2\pi}e^{i\theta}$, let us define
$f_{r,\theta}=f_\zeta$ by
\begin{eqnarray}
\widehat{f_\zeta}(\xi)&=&\left(1-\frac{\xi}{\zeta}\right)\left(1+\frac{\xi}{\zeta}\right)
\left(1-\frac{\xi}{\overline{\zeta}}\right)\left(1+\frac{\xi}{\overline{\zeta}}\right)\widehat{\ffi}(\xi)
\notag\\
&=&\left(1+\frac{2\cos 2\theta}{r^2}(2i\pi \xi)^2+\frac{1}{r^4}(2i\pi \xi)^4\right)\widehat{\ffi}(\xi).
\label{eq:struct}
\end{eqnarray}
Then $f_\zeta$ is positive definite.

Further, as $\ffi$ is smooth, and as $\widehat{\partial^k\ffi}(\xi)=(2i\pi \xi)^k\widehat{\ffi} (\xi)$
we get that
$$
f_\zeta(x)=\left(1+\frac{2\cos 2\theta}{r^2}\partial^2+\frac{1}{r^4}\partial^4\right)\ffi(x).
$$
The computations are easily justified by the fact that $\ffi$ is of class $\cc^4$.
Note that, for $|x|\leq 2$,
$$
\ffi''(x)=\frac{4}{105}(3x^4+18|x|^3+30|x|^2-12|x|-8)(2-x)^3
$$
and
$$
\ffi^{(4)}(x)=\frac{8}{5}(3x^4+6|x|^3-8x^2-16|x|+8)(2-x).
$$

We will now take $\theta=\pi/4$. Actually, we believe that for
every $0<\theta<\pi/2$ there is a unique $r$ such that $f_\zeta$
is extremal.

We then have
$$
f_\zeta(x)=\ffi(x)+\ffi^{(4)}(x)/r^4.
$$
It is not hard to see that $\ffi$ is decreasing and positive on
$(0,2)$ and that $\ffi^{(4)}$ is decreasing on
$(0,\sqrt{2-2/\sqrt{3}})$ and on $(\sqrt{2+2/\sqrt{3}},2)$
increasing on $(\sqrt{2-2/\sqrt{3}},\sqrt{2+2/\sqrt{3}})$,
positive on $(0,x_1)\cup(x_2,2)$ and negative on $(x_1,x_2)$ where
$x_1=0.441...$ and $x_2=1.462...$. In particular, for $r$ big
enough ($r=4$ will do), $f_{r,\pi/4}$ is positive on $(0,2)$ and
for $r$ small enough ($r=3$ will do), $f_{r,\pi/4}$ has two zeroes
on $(0,2)$. Therefore, there exists a unique $r$ such that
$f_{r,\pi/4}$ has exactly one double zero on $(0,2)$. A numerical
computation shows that $r\simeq 3.342775$. Let us denote this zero
by $x_\zeta$. A computer computation shows that $x_\zeta\simeq
1.303$.

We will need a bit more information. Let $\xi=\frac{\rho}{2\pi}e^{i\psi}$. Assume that $f_\xi\geq0$ and that there exists
$x_\xi\in (0,2)$ such that $f_\xi(x_\xi)=0$. Then, as $f_\xi$ is a polynomial on $(0,2)$,
we also have $f_\xi^\prime(x_\xi)=0$ that is
\begin{equation}
\left\{\begin{matrix}
\ffi(x_\xi)+2\frac{\cos 2\psi}{\rho^2}\ffi''(x_\xi)+\frac{1}{\rho^4}\ffi^{(4)}(x_\xi)&=&0\\
\ffi'(x_\xi)+2\frac{\cos 2\psi}{\rho^2}\ffi^{(3)}(x_\xi)+\frac{1}{\rho^4}\ffi^{(5)}(x_\xi)&=&0
\end{matrix}\right..
\label{eq:syst}
\end{equation}
A computer plot shows that $\ffi''\ffi^{(5)}-\ffi^{(3)}\ffi^{(4)}\not=0$ on $[0,2)$ so that,
taking the appropriate linear combination of both equations, we get
$$
2\frac{\cos 2\psi}{\rho^2}=\frac{\ffi'(x_\xi)\ffi^{(4)}(x_\xi)-\ffi(x_\xi)\ffi^{(5)}(x_\xi)}{
\ffi''(x_\xi)\ffi^{(5)}(x_\xi)-\ffi^{(3)}(x_\xi)\ffi^{(4)}(x_\xi)}.
$$
and
$$
\frac{1}{\rho^4}=\frac{\ffi'(x_\xi)\ffi''(x_\xi)-\ffi(x_\xi)\ffi^{(3)}(x_\xi)}
{\ffi^{(4)}(x_\xi)\ffi^{(3)}(x_\xi)-\ffi^{(5)}(x_\xi)\ffi''(x_\xi)}.
$$

In particular, $\rho$ and $\psi$ are uniquely determined by the
point where $f_\xi$ and its derivative vanish.

Let us now write $f_\zeta=g+h$ with $g,h$ positive positive
definite (either $L^2$ or continuous). Then, from Theorem
\ref{th.compextray}, $g$ and $h$ have at most $4$ complex zeroes
and, from the proof of that theorem and the remark following it,
we know that $g=\lambda f_{\xi}$ for some $\xi=\rho
e^{i\psi}\in\C$ and $\lambda>0$. But $0\leq g=\lambda f_\xi\leq
f_\zeta$ implies that $g$ must also have a double zero at
$x_\zeta$. Hence, by the previous argument, we must have
$\xi=\zeta$, and $g=\lambda f$.

\end{proof}

\section{Hermite functions}
In this section we restrict our attention to the one-dimensional situation.

\subsection{Preliminaries}\ 

Hermite functions are functions of the form $P(x)e^{-\lambda x^2}$
with $P$ a polynomial and $\lambda>0$. They satisfy many enjoyable
properties, in particular, they provide the optimum in many
\emph{uncertainty principles}, that is, their time-frequency
localization is optimal ({\it see e.g.} \cite{FS,HJ} and the
references therein).

Let us define the \emph{Hermite basis functions} by
$$
h_k(x)=\frac{2^{1/4}}{\sqrt{k!(4\pi)^{k}}}\,e^{\pi x^2}\partial^k
e^{-2\pi x^2}, \qquad k=0,1,2,\ldots
$$
It is well-known that $(h_k)_{k=0,1,\ldots}$  form an orthonormal
basis of $L^2(\R)$, and  $h_k(x)=c_kH_k(x)e^{-\pi x^2}$ with $H_k$
a real polynomial of degree $k$ with highest order term
$(2\sqrt{\pi}x)^k$ and $c_k$ a normalization constant that is not
relevant here. A simple computation shows that
\begin{equation}
\label{eq:hermite}
H_0(x)=1,\ H_2(x)=4\pi x^2-1\mbox{ and }
H_4(x)=(4\pi x^2)^2-6(4\pi x^2)+3.
\end{equation}

Finally, let us recall that the Hermite basis functions are
eigenvectors of the Fourier transform: $\widehat{h_k}=(-i)^kh_k$.
It immediately results that if the Hermite function
$H(x)=P(x)e^{-\pi x^2}$ is positive positive definite, then the
degree of $P$ is a multiple of $4$. Indeed, $H$ has to be even, so
that $P$ is real even. Then, if we expand $H$ in the Hermite
basis, we get that $\dst H(x)=\sum_{j=0}^k\alpha_j
H_{2j}(x)e^{-\pi x^2}$ where the $\alpha_j$'s are real and
$\alpha_k\not=0$. The Fourier transform of $H$ is then given by
$\dst\widehat{H}(\xi)=\sum_{j=0}^k(-1)^j\alpha_j
H_{2j}(\xi)e^{-\pi \xi^2}$. But, the highest order terms of the polynomial factor of $H$
and $\widehat{H}$ are then respectively $\alpha_k(4\pi x^2)^k$ and
$(-1)^k\alpha_k(4\pi \xi^2)^k$. Checking that $H$ and
$\widehat{H}$ stay non-negative when $x$ and $\xi$ go to infinity
suffices to see that $k$ is even.

\subsection{Gaussians are extremals}\ \\
Let us now turn to properties of extremals among Hermite functions. The first
result is a simple consequence of Hardy's Uncertainty Principle.

\begin{proposition}\label{prop:inthermite}\ \\
Let $\lambda>0$ and $P$ be a polynomial of degree $N$. Assume that
$f(x)=P(x)e^{-\lambda\pi x^2}$ is a positive positive definite
Hermite function. Then
$$
I(f)\subset\{Q(x)e^{-\lambda \pi x^2},\ Q \ \mbox{a polynomial of
degree }\leq N\}.
$$
In particular, if $f(x)=e^{-\lambda\pi  x^2}$, then $f$ is an
extremal ray generator in $\Omega^+$.
\end{proposition}

The second part of this proposition (and its proof) seems
well-known, {\it see} \cite{Bor}. The proof given here is only a
slight improvement.

\begin{proof} The second part of the proposition immediately follows from the first one.

Let $P$ be a polynomial of degree $N$ such that
$f(x)=P(x)e^{-\lambda\pi  x^2}\in\Omega^+$ and assume that
$P(x)e^{-\lambda\pi  x^2}=g(x)+h(x)$ with $g,h\in\Omega^+$. Then,
as $h(x)\geq 0$, $0\leq g(x)\leq P(x)e^{-\lambda\pi  x^2}$.
Further there exists a polynomial $\tilde{P}$ of degree $N$ such
that $\tilde{P}(\xi)e^{-\pi
\xi^2/\lambda}=\widehat{f}(\xi)=\widehat{g}(\xi)+\widehat{h}(\xi)$
so that, as $\widehat{h}\geq 0$, $0\leq \widehat{g}(\xi)\leq
\tilde{P}(\xi)e^{-\pi \xi^2/\lambda}$.

In particular, there exists a constant $C$ such that $|g(x)|\leq
C(1+|x|)^Ne^{-\lambda\pi  x^2}$ and $|\widehat{g}(\xi)|\leq
C(1+|\xi|)^Ne^{-\pi\xi^2/\lambda}$. From Hardy's Uncertainty
Principle (\cite{Ha} {\it see e.g.} \cite{FS,HJ} and \cite{BDJ,De1,De2} for
generalizations), there exists a polynomial $Q$ of degree at most
$N$ such that $g(x)=Q(x)e^{-\lambda\pi  x^2}$ and therefore
$h(x)=(P-Q)(x)e^{-\lambda\pi  x^2}$.
\end{proof}

Let us conclude with the following lemma:

\begin{lemma}
\label{lem:propextherm} Let $\lambda>0$ and $P$ be a polynomial of
degree $4N\geq 4$ and assume that $f(x)=P(x)e^{-\pi\lambda x^2}$
is a positive positive definite Hermite function. Write
$\widehat{f}(\xi)=\widetilde{P}(\xi)e^{-\pi\xi^2/\lambda}$.

--- If $f$ is an extremal ray generator, then either $P$ or $\widetilde{P}$ has at least 4 real zeroes.

--- If $P$ or $\widetilde{P}$ has $4N$ real zeroes, then $f$ is an extremal ray generator.
\end{lemma}

\begin{proof} It is enough to prove the lemma with $\lambda=1$.
Assume that $P$ and $\widetilde{P}$ are non-negative and have no
zeroes, then there exists a constant $C>0$ such that, $P(x)-C$ and
$\widetilde{P}(\xi)-C$ are still non-negative. It follows that
$\frac{P(x)-C}{2}e^{-\pi x^2}$ and its Fourier transform,
$\frac{\widetilde{P}(\xi)-C}{2}e^{-\pi \xi^2}$ are both positive.
As $f(x)=\frac{P(x)-C}{2}e^{-\pi x^2}+\frac{P(x)+C}{2}e^{-\pi
x^2}$ we conclude that $f$ is not an extremal ray generator. It
follows that $P$ or $\widetilde{P}$ has at least one zero on the real line. As they
are positive polynomials, any real zero must have even multiplicity. As
these polynomials are even and non-zero at the origin, we get that
$P$ or $\widetilde{P}$ indeed has at least $4$ real zeroes.

For the second part, it is enough to consider the case when $P$
has $4N$ real zeroes. If $f=g+h$ with $g,h\in I(f)$ then, from
Proposition \ref{prop:inthermite}, we get that $g(x)=Q(x)e^{-\pi
x^2}$ and $h(x)=R(x)e^{-\pi x^2}$ with $Q,R$ polynomials of degree
at most $4N$. As $h\geq 0$, $0\leq g\leq f$ thus $0\leq Q\leq P$.
It follows that a real zero of $P$ is also a real zero of $Q$. As
$Q$ has not higher degree than $P$, $Q=cP$ with $0\leq c\leq 1$
and then $g=cf$, $h=(1-c)f$.
\end{proof}

\subsection{Extremal ray generators among Hermite functions of higher degree}\

From the previous section, we know that if a Hermite function
$f(x)=P(x)e^{-\pi x^2}$ is positive positive definite then the
polynomial has degree a multiple of 4. Moreover, the interval
generated by $f$ consists of Hermite functions of not higher
degree. Thus, if there exists a polynomial $P$ of degree $4q$ such
that $P(x)e^{-\pi x^2}$ is positive positive definite, then there
exist extremals of the same degree. Indeed, we just have to
consider the finite dimensional cone of positive positive definite
Hermite functions of degree $4q$, which is then non-empty, and is
thus the positive span of its extremal rays.

We will now characterize all positive
positive definite Hermite functions of degree $4$ and the
extremals among them. It is enough to consider $f$ of the form
$f(x)=(H_0(x)+2aH_2(x)+bH_4(x))e^{-\pi x^2}$. Using the fact that
the Hermite basis consists of eigenvalues of the Fourier
transform, we get that
$\widehat{f}(\xi)=(H_0(\xi)-2aH_2(\xi)+bH_4(\xi))e^{-\pi \xi^2}$.

We thus aim at characterizing $a,b$ for which $P(x)=H_0(x)+2aH_2(x)+bH_4(x)\geq 0$.
But $H_0(x)+2aH_2(x)+bH_4(x)=1+2a+3b-8\pi(a+3b)x^2+16b\pi^2x^4$.
Setting $X=4\pi x^2$, we thus ask whether
$\widetilde{P}(X):=1+2a+3b-2(a+3b)X+bX^2\geq0$ for all $X\geq 0$.

The first condition is that $\widetilde{P}(0)\geq 0$ that is $1+2a+3b\geq 0$ and that
$\lim_{X\to+\infty}\widetilde{P}(X)\geq 0$, that is $b\geq 0$. Next, we want that
$\widetilde{P}$ has no single root in $]0,+\infty)$. Thus, either $a+3b\leq 0$ or
$(a+3b)^2-b(1+2a+3b)\leq 0$
which we may rewrite as
\begin{equation}
\label{eq:ellipse}
\Bigl(a+2b\Bigr)^2+2\Bigl(b-\frac{1}{4}\Bigr)^2\leq\frac{1}{8}.
\end{equation}
This is the equation of an ellipse $\mathcal{E}_+$,
which passes through the point $(0,0)$ where it is tangent
to the line of equation $b=0$ and through the point $(-1,1/3)$
where it is tangent to the line of equation $1+2a+3b=0$.

In this form, we see that the set $\mathcal{D}_+$ of all $(a,b)$'s
for which $P\geq 0$ is the union of the ellipse of equation
\eqref{eq:ellipse} and the triangle formed by the lines $b=0$,
$a+3b=0$ and $1+2a+3b=0$, that is the triangle with vertices
$(0,0)$, $(-1,1/3)$ $(-1/2,0)$.

If we want $f$ to be positive positive definite, then both $(a,b)$
and $(-a,b)$ have to belong to $\mathcal{D}_+$. It is easy to see
that this reduces to the intersection of the two ellipses. More precisely we get\,:

\begin{proposition}
\label{prop:hermite4}
The set of all Hermite functions of order $4$ that are positive positive definite is given by
$f_{a,b}:=(H_0+2aH_2+bH_4)e^{-\pi x^2}$ where $(a,b)$ belongs to the region $\mathcal{D}$
parametrized by the following\,:

--- either $a\geq 0$ and  $(a+2b)^2+2(b-1/4)^2\leq 1/8$.

--- or $a\leq 0$ and $(-a+2b)^2+2(b-1/4)^2\leq 1/8$.

Moreover $f_{a,b}$ is an extremal ray generator of $\Omega^+$ if
and only if $(a,b)\in\partial\mathcal{D}$, where
$\partial\mathcal{D}$ is the boundary of $\mathcal{D}$.
\end{proposition}

\begin{figure}[ht]
\begin{center}
\scalebox{0.6}{\includegraphics[height=5cm,width=12cm]{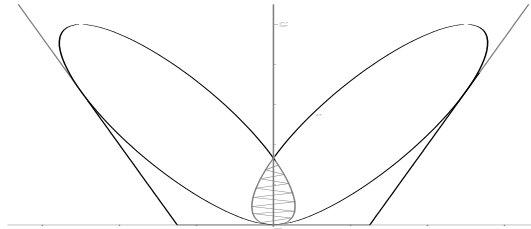}}
\end{center}
\label{figident}
\caption{The set of all $(a,b)$'s for which $f_{a,b}$ is positive positive definite.}
\end{figure}

The last part results directly from the fact that $\mathcal{D}$ is
convex and that all its boundary points are points of curvature
(excepted two) and are thus extreme points of $\mathcal{D}$. The
fact that $f_{a,b}$ is then an extremal of $\Omega^+$ is a direct
consequence of the previous discussion.

It should also be noted that if $f_{a,b}$ is an extremal ray generator, then

--- either $f_{a,b}$ has only real zeroes and $\widehat{f_{a,b}}$ is positive,
in which case we will say that $f_{a,b}$ is of the \emph{time type}

--- or $\widehat{f_{a,b}}$ has only real zeroes and $f_{a,b}$ is positive,
in which case we will say that $f_{a,b}$ is of the \emph{frequency type}.

\begin{corollary}
Assume that $f_{a_1,b_1},\ldots,f_{a_N,b_N}$ are all of the time type,
then $f=\dst\prod_{i=1}^Nf_{a_i,b_i}$ is an extreme ray generator in $\cc$.

If $f_{a_1,b_1},\ldots,f_{a_N,b_N}$ are all of the frequency type,
then  $f=\dst\star_{i=1}^Nf_{a_i,b_i}$ is an extreme ray generator in $\cc$.
\end{corollary}


\begin{proof}
Let us first assume that all the $f_{a_i,b_i}$'s are of time type
and consider $f=\dst\prod_{i=1}^Nf_{a_i,b_i}$. If this is
extremal, then $\widehat{f}=\star_{i=1}^Nf_{-a_i,b_i}$ is also
extremal.

Let us write $f=P\gamma^N=g+h$ with $P$ a polynomial, $\gamma$ the
standard gaussian, $g,h\in\Omega^+$. From Proposition
\ref{prop:inthermite}, we know that both $g$ and $h$ are Hermite
functions, $g=G\gamma^N$ and $h=H\gamma^N$ with $G,H$ positive
polynomials of degree $\leq 4N$.

But, as $h\geq 0$, $0\leq G\leq P$ and every real zero of $P$ is thus a real zero
of $G$. As $P$ is assumed to have $4N$ real zeroes, $G$ has $4N$ real zeroes.
As $G$ is of degree at most $4N$, $G=\lambda P$ thus $g=\lambda f$ and $h=(1-\lambda)f$.

If all the $f_{a_i,b_i}$'s are of frequency type, the same
argument applies on the Fourier side.
\end{proof}

\section{Some important classes of positive positive definite functions}
\label{sec:criteria}

In this section, we will show that several criteria available in the literature
are particular cases of Choquet representation of the cone $\Omega_d^+$.
We will appeal several times to lemma \ref{lem:extcrit}.

\subsection{Non-negative convex functions}
\label{sec:polya}\ \\
A well-known class of positive definite functions is that of
\emph{P\'olya type} functions ({\it see} \cite{Po}). More
precisely, $f$ is said to be of P\'olya type if $f$ is even,
continuous, convex on $[0,+\infty)$ and $f(x)\to0$ when $x\to
\infty$. It is well-known that $f$ is of P\'olya type if and only
if there exists a positive bounded measure $\nu$ on $[0,+\infty)$
such that, for all $x\in\R$.
\begin{equation}
\label{eq:polya}
f(x)=\int_0^{+\infty}(1-|x|/t)_+\,\mbox{d}\nu(t).
\end{equation}
Note that Lebesgue's Theorem ensures that a function $f$ given by \eqref{eq:polya}
is continuous as soon as the measure $\nu$ appearing in that formula is bounded.

Further, one may show that $f$ is a P\'olya function then $f$ is a
Fourier transform of a measure of the form $p(\xi)\,\mbox{d}\xi$
with $p$ continuous on $\R\setminus\{0\}$. More precisely, $\nu$
and $p$ are related by
$$
p(\xi)=c\int_0^{+\infty}\left(\frac{\sin(\pi t\xi)}{t\xi}\right)^2t\,\mbox{d}\nu(t)
$$
where $c$ is some constant (and $\frac{\sin 0}{0}=1$).

Note that everything is straightforward when $f$ has sufficient
smoothness and decay, using integrations by part. For a simple
proof in the general case and further references, we refer to
\cite{Sa}.

Finally, an easy computation shows that $(1-|x|/t)_+=c\chi_{[-t/2,t/2]}*\chi_{[-t/2,t/2]}(x)$
which is an extremal ray generator according to Theorem \ref{th.compextray}.
We may thus rewrite \eqref{eq:polya} in the form
\begin{equation}
\label{eq:convconv}
f(x)=\int_0^{+\infty} \chi_{[-t/2,t/2]}*\chi_{[-t/2,t/2]}(x)\,\mbox{d}\nu(t).
\end{equation}
In this form, we immediately see that we are in the situation of Lemma \ref{lem:extcrit}.

\begin{proposition}
\label{prop:conv}\ \\
Let $f$ be a P\'olya type function. Then $f$ is positive positive
definite, and it is an extremal ray generator if and only if
$f=\chi_{[-r/2,r/2]}*\chi_{[-r/2,r/2]}$ for some $r>0$.
\end{proposition}

\begin{example}\ \\
--- $\psi(t)=\frac{\ln(e+|t|)}{1+|t|}$ and $\psi(t)=\frac{1}{(1+|t|)^\alpha}$, $\alpha>0$ are of
P\'olya type and are thus not extremal ray generators.

\noindent--- It is easy to show that $e^{-|x|^p}$ satisfies the
hypothesis of the theorem for $0<p\leq 1$, in particular, they are
not extremal ray generators. We will show below that this stays
true for $1<p<2$ as well.
\end{example}

\subsection{Generalizations by Gneiting}
\label{sec:gneiting}\ \\
P\'olya's criterion has been extended recently by Gneiting
\cite{Gn1,Gnproc} and by C. Hainzl and R. Seiringer \cite{HS}. For
instance, \cite{Gn1} characterizes those functions which can be
written in the form
\begin{equation}
\label{gneiting}
\psi(x)=\int_0^{+\infty}w(xt)\,\mbox{d}\mu(t)
\end{equation}
where $w(x):=(1-x^2)_+*(1-x^2)_+$ and $\mu$ is a positive bounded
measure on $[0,+\infty)$. Note that $w$ is the positive positive
definite function introduced by Wu in the study of radial basis
function interpolation and that
$w(x)=c(1-|x|)_+^3(1+3/2|x|+x^2/4)$ where $c$ is a constant.
This function also plays an important role in the work of Wendland \cite{We}
on positive definite functions with compact support of optimal smoothness.

More precisely, Gneiting showed that a function
$\psi(t)=\ffi(|t|)$ can be written in the form \eqref{gneiting} if
and only if $\ffi$ satisfies
\begin{enumerate}
\renewcommand{\theenumi}{\roman{enumi}}
\item $\ffi$ is twice continuously differentiable,

\item $\ffi(0)>0$ and $\ffi(x)\to0$ when $x\to+\infty$,

\item $\dst\frac{1}{t}\bigl(\sqrt{t}\ffi''(\sqrt{t})-\ffi'(\sqrt{t})\bigr)$ is convex.
\end{enumerate}
Note that such a function is necessarily non-negative.

The same proof as in the previous section then shows the following:

\begin{proposition}\ \\
A function $\psi$ of the form \eqref{gneiting} (or equivalently a
function $\psi$ given by $\psi(t)=\ffi(|t|)$ with $\ffi$
satisfying $i)-iii)$ above) is positive positive definite, and it
is an extremal ray generator in the set of positive positive
definite functions if and only if there exists $c>0$ such that
$\psi(t)=w(ct)$.
\end{proposition}

\begin{example}\ \\
As already noted by Gneiting, the following functions satisfy \eqref{gneiting} and are therefore
not extremal ray generators:
$$
i)\ \psi(t)=\frac{1}{1+|t|^\beta},\ 0<\beta<1.877...
,\quad ii)\ \psi(t)=(1+\gamma |t|+t^2) \exp(-|t|),\ 0<\gamma<1/4
$$
$$
iii)\ \psi(t)=(1-|t|)_+^3(1+3|t|)
\quad iv)\  \psi(t)=(1-|t|^\lambda)_+^3
$$
The first one is known as Linnik's function and is of P\'olya type
for $\beta<1$. The third one has been introduced by Wendland and
has interesting optimal smoothness properties. The last one is
called Kuttner's function and seems to be Gneiting's original
motivation in the above characterization.

Finally, $\psi(t)=\exp(-|t|^\beta)$ satisfies the hypothesis of the theorem
for $\beta<1.84170$ thus improving the domain of non-extremality found in the previous section.
\end{example}

\subsection{Bernstein functions}\ \\
Let us recall that a non-negative function $g$ on $[0,+\infty)$
is called \emph{completely monotonic} if it is infinitely
differentiable on $(0,+\infty)$ and, for all $k\in\N$ and
$x\in(0,+\infty)$,
$$
(-1)^kg^{(k)}(x)\geq 0.
$$
Completely monotonic functions have remarkable applications in
different branches. For instance, they play a role in potential
theory \cite{BF}, probability theory \cite{Bon,Fe,Ki}, physics
\cite{Da}, numerical and asymptotic analysis \cite{Fr,Wim}, and
combinatorics \cite{Ba}. A detailed collection of the most
important properties of completely monotonic functions can be
found in \cite[Chapter IV]{Wid}, and in an abstract setting in
\cite{BCR}.

The celebrated theorem of Bernstein ({\it see} \cite[Chapter III
Setion 2]{Be}, \cite[Chapter 18, Section 4]{Fe} or \cite[page
161]{Wid}) states that every  completely monotonic function which
is continuous at zero is the Laplace transform of a positive bounded measure
on $[0,+\infty)$. In other words, $g$ is completely monotonic and
continuous at $0$ if and only if there exists a (necessarily
unique) positive bounded measure $\mu$ on
$[0,+\infty)$ such that, for every $x\geq 0$,
\begin{equation}
\label{eq:bernstein}
g(x)=\int_0^{+\infty}e^{-tx}\,\mbox{d}\mu(t).
\end{equation}
A particular case of Lemma \ref{lem:extcrit}, of which the first
part is well-known, is then the following:

\begin{proposition}\ \\
Let $g$ be a non-negative, completely monotonic function, continuous at $0$.
Let $f$ be defined on $\R^d$ by $f(x)=g(|x|^2)$. Then $f$ is a positive positive definite function
and is an extremal ray generator if and only if $f$ is a Gaussian.
\end{proposition}

\begin{example}\ \\
--- Let $\lambda>0$ and $0<\alpha\leq 1$.
It is well-known and easy to show that $g$ defined by
$g(t)=e^{-\lambda t^\alpha}$ is completely monotonic. It follows
that $e^{-\lambda|x|^p}$ is positive positive definite for
$0<p\leq 2$ and an extremal ray generator only for $p=2$.

--- Let $\alpha\not=0$ and $\beta>0$. It is easy to show that $g(t)=(t+\alpha^2)^{-\beta}$
(the so-called \emph{inverse multiquadrics}) is completely
monotonic. It follows that $(x^2+\alpha^2)^{-\beta}$ is positive
positive definite but not an extremal ray generator.

--- Recently, the so-called Dagum family $D_{\beta,\gamma}(x)=1-\left(\frac{x^\beta}{1+x^\beta}\right)^\gamma$
has been introduced in \cite{PMZP,MPN} where it was shown that the
Dagum class allows for treating independently the fractal
dimension and the Hurst effect of the associated weakly stationary
Gaussian RF, by using the procedure suggested in \cite{GS}. In
\cite{BMP}, the authors further showed that for certain
parameters, this function is completely monotonic.
\end{example}

Note that completely monotonic functions allow to construct radial
positive definite functions in \emph{any} dimension. There is a
converse to this. More precisely, let $f_0$ be a continuous
function on $[0,+\infty)$ and define $f_d$ on $\R^d$ by
$f_d(x)=f_0(|x|)$. Observe that if $f_d$ is positive definite then
$f_{d'}$ is also positive definite for any $d'<d$. A famous
theorem of Schoenberg \cite{Sc} (see also \cite{SvP} for a more
modern proof and further references or \cite{Bax} for another
proof) states that if $f_d$ is positive definite for any $d$, then
$f_0(x)=g(x^2)$ with $g$ completely monotonic.

\section{Conclusion}

In this paper we have investigated the extremal ray generators of
the cone of continuous positive positive definite functions for which
we have described two important classes of such extremals. 

At this stage, we would first like to stress that, up to minor modifications, most of our results
stay true for positive positive definite $L^2$ functions or even tempered distributions.
Let us recall that in this case, one can not define positive definiteness via \eqref{def.intposdef}.
One thus replaces this condition by one that is equivalent for continuous functions.
To do so, notice that if $f$ is continuous, then $f$ is positive definite if and only if, for
every smooth function $\Phi\in\ss(\R^d)$ in the Schwartz-class,
\begin{equation}
\label{def.intposdef}
\int_{\R^d}\int_{\R^d} f(x-y)\Phi(x)\overline{\Phi(y)}\,\mbox{d}x\,\mbox{d}y\geq 0.
\end{equation}
Note that this makes sense if $f$ is only assumed to be in
$L^\infty$, but in that case it is well-known ({\it e.g.} \cite{Be,Ch,Ph})
that  implies that $f$ is almost everywhere
equal to a continuous function. Further, we may extend
\eqref{def.intposdef} to $T\in \ss'$ and take this to be the
definition of a positive definite tempered distribution\,: a
tempered distribution $T$ is positive definite if for every
$\Phi\in\ss(\R^d)$, $\scal{T,\Phi*\Phi^*}\geq 0$ where
$\Phi^*(x)=\overline{\Phi(-x)}$. Schwartz \cite{Sch} extended
Bochner's theorem by characterizing positive definite tempered
distributions as the Fourier transforms of positive tempered
measures, that is
$\scal{T,\Phi}=\int_{\R^d}\widehat{\Phi}(\xi)\,\mbox{d}\mu(\xi)$
where $\mu$ is a positive measure such that, for some $\alpha>0$,
$\int_{\R^d} (1+|\xi^2|)^{-\alpha}\,\mbox{d}\mu(\xi)<+\infty$. A
proof based on Choquet Theory may be found in \cite{To2}. As a
particular case, one gets that if $T\in L^2$, then $T$ is positive
definite if and only if $\widehat{T}\geq0$. For further extensions
of the notion of positive definiteness, we refer to \cite{Ste}.

\medskip

Unfortunately, several open problems remain, but some features seem to have become clear:

--- A full classification of the extremals is probably impossible and no reasonable
conjecture can be stated at this stage. For instance, the
existence of an extremal compactly supported function with Fourier
transform having complex zeroes leaves little hope for a full
description of compactly supported extremals. As for the extremals
inside the Hermite class, we have a full description only up to
degree $4$. For higher degree, even though we have been able to
construct extremals, we are far from a full description. The main
reason for this, is that we are unable to decide the answer of the
following question:

\medskip

\noindent{\bf Question 1.}\ \\
{\sl a) Is the product of two extremals an extremal?\\
b) What about the convolution of two extremals (if it makes
sense)?}

\medskip

We have no hint of what the answer might be and will therefore not
propose a conjecture.

\medskip

--- There are many criteria allowing to decide whether a function is
positive definite. As we have seen, these criteria actually
characterize the functions that are mixtures of the scales of a
single extremal inside the class of positive positive definite
functions. We are thus tempted to ask for more criteria, for
instance:

\medskip

\noindent{\bf Question 2.}\ \\
{\sl a) Find a criterion that allows to decide whether a function
is the mixture of the scales
of $m_\alpha*m_\alpha$, where $\alpha$ may, or may not be fixed.\\
b) Find a criterion that allows to decide whether a function is
the mixture of extremals
of scales of extremal Hermite functions of order 4.\\
c) Find a criterion for positive definiteness that is not given in
terms of a mixture of positive positive definite functions, thus
also characterizing positive definite functions that are not
necessarily non-negative.}

\medskip

--- Finally, all the extremals we found, except the gaussians,
either have at least a zero, or their Fourier transform has at
least a zero. Borisov \cite{Bor} conjectured that the only
extremals in $L^2(\R)\cap\cc(\R)$ that are strictly positive and
that have strictly positive Fourier transform are the gaussian
functions $e^{-a\pi x^2}$, $a>0$. We are inclined to believe that
this is false. This is linked to Question 1 above which leads us
to the following more precise question:

\medskip

{\bf Question 3.}\ \\
{\sl Consider the two functions
$f=\chi_{[-1/2,1/2]}*\chi_{[-1/2,1/2]}$ and $\gamma=e^{-\pi x^2}$,
and recall that they are extremal positive positive definite
functions in $\Omega_1^+$. Is $f*\gamma$ extremal in $\Omega_1^+$
and, if so, is $\gamma (f*\gamma)$ extremal in $\Omega_1^+$?}

\medskip

A positive answer to the second question would of course provide a
counterexample to Borisov's conjecture.

\medskip

Overall, we hope this paper will help reviving P. Levy's request
to study the cone of positive positive definite functions.

\end{document}